\documentstyle[amsmath,amssymb,11pt]{article}

\newtheorem{theorem}{\bf Theorem}[section]

\newtheorem{lemma}{\bf Lemma}

\begin{document}

{\rm

\vspace{1cm}

\begin{center}

{\large  Inversion Theorem for Bilinear Hilbert Transform}

\vspace{0.5cm}

Aneta Bu\v ckovska,\footnote{University of Skopje, R. Macedonia} Stevan Pilipovi\'
c\footnote{University of Novi Sad,
 Serbia}
Mirjana Vukovi\'
c\footnote{University of Sarajevo, Bosnia and Herzegovina}

\vspace{0.5cm}

\end{center}

\begin{abstract}
  An approximation  result for the bilinear Hilbert transform
is proved and used for the inversion of the bilinear Hilbert transform.
Also, p-Lebesgue points $(p\geq 1)$ are analyzed.
\end{abstract}

\vspace{0.2cm}

\vspace{0.2cm}


\vspace{0.5cm}

\section{Introduction}

The aim of this paper is to give an  inversion theorem for the
bilinear Hilbert transform (BHT) defined  in appropriate classes
of functions and distributions. More precisely, for functions, the
product $f(x)g(x)$ is obtained as the inversion of the BHT at
Lebesgue points (Theorem 1.1).

 In  several papers Lacey and Thiele ([4]-[6]) had studied the
continuity of  BHT
$$
H_{\alpha}(f,g)(x)=\lim_{\varepsilon\rightarrow 0}{\int_{0<\varepsilon\leq|t|} f(x-y)g(x+\alpha y){dy\over y}},\quad
\alpha \in{\mathbb{R}}\setminus \{0,-1\}\;,$$ where $f \in
L^{2}(\mathbb{R})$ and $g\in L^\infty(\mathbb{R}),$ respectively
$f \in L^{p_1}(R)$ and $g\in L^{p_2}(\mathbb{R}),$ $1<p_1,
p_2<\infty. $ Their main result is the affirmative answer on the
Calderon conjecture, first for $p_1=2, p_2= \infty$ ([5]), then
for $p_1$, $p_2\in (1,\infty)$. Let $2/3<p=\frac{p_1p_2}{p_1+p_2}$
or $p_1=2,\;p_2=\infty$ and $p=2$.  Then their main result is $
 ||H_\alpha(f,a)||_{L^p} \leq C ||f||_{L^{p_1}}||a||_{L^{p_2}},
 \; f\in L^{p_1},\;
 a \in L^{p_2},
$ where $C>0$ depends on $\alpha, p_1, p_2$. We refer to [7] and
the references therein
 for further reading on multi-linear operators given by singular multipliers.

The bilinear Hilbert transform $H_{\alpha}:L^{2}\times L^\infty
\rightarrow L^2 $ respectively, $H_{\alpha}:L^{p_1}\times L^{p_2}
\rightarrow L^{p}, $ was extended in  [1] to ${\cal
D}'_{L^{2}}\times {\cal D}_{L^{\infty}} \rightarrow {\cal
D}'_{L^2},$ respectively, ${\cal D}'_{L^{q}}\times {\cal
D}_{L^{p_2}} \rightarrow {\cal D}'_{q_1},$ (with suitable
parameters) as a  hypocontinuous, respectively, continuous
mapping. The trilinear Hilbert transform on ${\cal D}_{L^p}\times
{\cal D}_{L^q}\times {\cal D}_A\rightarrow {\cal D}_{L^r}$ is
studied in [3]. This analysis is based on [7]. Moreover the BHT of
ultradistributions is analyzed in [2].

Recall [9]: Let $f\in L^p_{loc}$, $p\geq 1$. Then $x\in
\mathbb{R}$ is a p-Lebesgue point of $f$, $x\in A^p_f$, if
$$\frac {1}{r}\int\limits_{|t|<r} {\vert f(x-t) - f(x)\vert}^p dt \rightarrow 0,{\mbox as }\;r\rightarrow 0.$$
If $p=1$, we will use notation $A_f^1=A_f$.

Our aim is to prove the following theorem.

\begin{theorem}
Let $f\in L^{2}(\mathbb{R})$, $g\in L^{\infty}(\mathbb{R})$,
respectively $f\in L^{p_1}(\mathbb{R})$, $g\in
L^{p_2}(\mathbb{R})$, ${\frac{1}{p_1}}+{\frac{1}{p_2}}\leq 1$. Let
$x\in A^2_f\bigcap A^\infty_g$, respectively, $x\in
A^{p_1}_f\bigcap A^{p_2}_g$. Then
$$
f(x)g(x)=\frac{i}{\pi}(\lim_{\varepsilon\rightarrow
0}H_{\alpha,\varepsilon}(f,g)(x)-H_{\alpha}(f,g)(x)),
$$
where
$$
H_{\alpha,\varepsilon}(f,g)(x)=\int\limits_{\mathbb{R}} f(x-t)g(x+\alpha
t){\frac{dt}{t+i\varepsilon}},\quad x\in {\mathbb{R}},\;\varepsilon\in(0,1).
$$
\end{theorem}
 This theorem is stated in [1].

{\bf Outline of the proof.}

Since
$$
\lim_{\varepsilon\rightarrow 0}H_{\alpha,\varepsilon}(f,g)(x)-H_{\alpha}(f,g)(x)=
$$
\begin{equation}
\label{a} \lim_{\varepsilon\rightarrow
0}\big(\int\limits_{-\infty}^{\infty} f(x-t)g(x+\alpha t)\frac
{t}{t^2+\varepsilon^2}dt- \int_{0<\varepsilon\leq|t|}
f(x-y)g(x+\alpha y){dy\over y}\big)
\end{equation}
$$
+
\lim_{\varepsilon\rightarrow 0}
\int_{0<\varepsilon\leq|t|} f(x-y)g(x+\alpha y){dy\over y}
-H_{\alpha}(f,g)(x)
$$
\begin{equation}
\label{b}
-i\lim_{\varepsilon\rightarrow 0}\int\limits_{-\infty}^{\infty}[f(x-t)g(x+\alpha
t)-f(x)g(x)]\frac {\varepsilon dt}{t^2+\varepsilon^2}
\end{equation}
$$
-i\pi
f(x)g(x),
$$
 we have to prove that (\ref{a}) and (\ref{b}) tend to zero as
$\varepsilon\rightarrow 0$, if  $x\in A^2_f\bigcap A^\infty_g$,
respectively, $x\in A^{p_1}_f\bigcap A^{p_2}_g.$

\vspace{05cm}

 For this proof  we need an appropriate
analysis of p-Lebesgue points. Section 2 is devoted to p-Lebesgue
points, Section 3 contains a preparation of Theorem 1.25, in [8], Chapter
I, in the context of BHT and finally, at the end of Section
3 Theorem 1.1 is proved. This theorem is used in the appendix for the  the inversion theorem
of
BHT in the spaces of distributions and ultradistributions.

\section{On p-Lebesgue points }

\begin{lemma} Let $1\leq p_2\leq p_1$ and $x\in A^{p_1}_f$. Then $x\in A^{p_2}_f$.
\end{lemma}
{\bf Proof:} Let $p=\frac{p_{1}}{p_2}$ and $q=\frac{p_1}{p_1-p_2}$.
Since $\frac{1}{p}+\frac{1}{q}=1$ we have
\begin{eqnarray}
\frac {1}{r}\textstyle\int\limits_{|t|<r} {\vert f(x-t) - f(x)\vert}^{p_2}dt & \leq &
\frac {1}{r^{{\frac{1}{p}}+ {\frac{1}{q}}}}\Bigl(\textstyle\int\limits_{|t|<r}{\vert
f(x-t) - f(x)\vert}^{p_1}dt\Bigr)^{\frac{1}{p}}\nonumber \\
&   & \cdot\Bigl(\textstyle\int\limits_{|t|<r}dt\Bigr)^{\frac{1}{q}}\nonumber \\
& = & \Bigl({\frac{1}{r}}\textstyle\int\limits_{|t|<r}{\vert f(x-t) - f(x)\vert}^{p_1}dt\Bigr)^{\frac{1}{p}}\nonumber
\end{eqnarray}
and this implies the assertion.
\begin{lemma} If $f\in L^{p_1}_{loc}$, $g\in L^{p_2}_{loc}$, where ${\frac{1}{p_1}}+{\frac{1}{p_2}}=1$ and
$x\in A^{p_1}_f \bigcap A^{p_2}_g,$ then $x\in A_{fg}$.
\end{lemma}
{\bf Proof:}
We have
\begin{eqnarray}
K & = & {\frac {1}{r}}\textstyle\int\limits_{|t|<r} {\vert f(x-t)g(x-t) - f(x)g(x)\vert}dt \nonumber \\
& \leq & {\frac{1}{r}}\textstyle\int\limits_{|t|<r}{\vert f(x-t) - f(x)\vert}\cdot \vert g(x-t)\vert dt \nonumber \\
&  & +\frac{1}{r}\textstyle\int\limits_{|t|<r}{\vert g(x-t) -
g(x)\vert\cdot\vert f(x)\vert}dt \nonumber \\
& = & I+J. \nonumber
\end{eqnarray}
For the first integral $I$ we have:
\begin{eqnarray}
I& \leq & \textstyle{\frac {1}{r}}\Bigl(\textstyle\int\limits_{|t|<r}{\vert f(x-t) -
 f(x)\vert^{p_1}dt\Bigr)^{\frac{1}{p_1}}
\Bigl(\textstyle\int\limits_{|t|<r}\vert g(x-t)-g(x)+g(x)\vert^{p_2}
}dt\Bigr)^{\frac{1}{p_2}} \nonumber \\
& \leq & \frac {1}{r^{{\frac{1}{p_1}}+{\frac{1}{p_2}}}}\Bigl(\textstyle\int\limits_{|t|<r}
{\vert f(x-t) - f(x)\vert}^{p_1}dt\Bigr)^
{\frac{1}{p_1}}\Bigl[\Bigl(\textstyle\int\limits_{|t|<r}\vert
g(x-t)-g(x)\vert^{p_2}dt\Bigr)^{\frac{1}{p_2}} \nonumber \\
&  & +\Bigl(\textstyle\int\limits_{|t|<r}{\vert g(x)\vert}^{p_2}dt\Bigr)^{\frac{1}{p_2}}\Bigr] \nonumber \\
& = & \Bigl(\textstyle{\frac{1}{r}}\textstyle\int\limits_{|t|<r}{\vert f(x-t) - f(x)\vert}^{p_1}dt\Bigr)^{\frac{1}{p_1}}
\Bigl[\Bigl({\frac{1}{r}}\textstyle\int\limits_{|t|<r}{\vert g(x-t) -
g(x)\vert}^{p_2}dt\Bigr)^{\frac{1}{p_2}} \nonumber \\
& & +\vert g(x)\vert\frac{{(2r)}^{\frac{1}{p_2}}}{r^{\frac{1}{p_2}}}\Bigr]. \nonumber
\end{eqnarray}
Because $x\in A^{p_1}_f \bigcap A^{p_2}_g$, it follows that
$${\frac{1}{r}}\int\limits_{|t|<r}{\vert f(x-t) - f(x)\vert}^{p_i}dt \rightarrow 0 \;{\mbox  as  }\; r\rightarrow 0,
\quad i=1,2.$$ So  integral $I\rightarrow 0$ as
$r\rightarrow 0$. Integral $J$ also tends to zero as $r\rightarrow
0$, according to Lemma 1 and we obtain that $K$ tends to zero and
that $x\in A_{fg}$.
\begin{lemma} Let $f\in L^{p_1}_{loc}$, $g\in L^{p_2}_{loc}$, where
${\frac{1}{p_1}}+{\frac{1}{p_2}}={\frac{1}{p_3}}<1$, $x\in
A^{p_1}_f\bigcap A^{p_2}_g$.  Then $x\in A^{p_3}_{fg}$.
\end{lemma}
{\bf Proof:} First we note that $fg\in L^{p_3}_{loc}$. We have
\begin{eqnarray}
K & = & {\frac {1}{r}}\int\limits_{|t|<r} {\vert f(x-t)g(x-t) - f(x)g(x)\vert}^{p_3}dt \nonumber \\
& \leq & {\frac {1}{r}}\int\limits_{|t|<r}{\vert g(x-t) - g(x)\vert^{p_3}\cdot \vert
f(x-t)\vert^{p_3}}dt \nonumber \\
&&+ \frac{1}{r}\int\limits_{|t|<r}{\vert f(x-t) - f(x)\vert^{p_3}\cdot\vert g(x)\vert^{p_3}}dt \nonumber \\
& = & I+J. \nonumber
\end{eqnarray}
Since ${\frac{p_3}{p_1}}+{\frac{p_3}{p_2}}=1,$ for the first
integral $I$ we have:
\begin{eqnarray}
I & \leq & {\frac {1}{r}}\Bigl(\int\limits_{|t|<r}{\vert f(x-t)\vert^{p_1}dt\Bigr)^{\frac{p_3}{p_1}}\cdot
\Bigl(\int\limits_{|t|<r}\vert g(x-t)-g(x)\vert^{p_2}
}dt\Bigr)^{\frac{p_3}{p_2}} \nonumber \\
& \leq &{\frac {1}{r}}\Bigl(\int\limits_{|t|<r}\vert f(x-t)-f(x)\vert^{p_1}dt+\int\limits_{|t|<r}\vert
f(x)\vert^{p_1}dt\Bigr)^{\frac{p_3}{p_1}}
\Bigl(\int\limits_{|t|<r}\vert g(x-t)-g(x)\vert^{p_2}dt\Bigr)^{\frac{p_3}{p_2}} \nonumber \\
& \leq & \frac {1}{r^{{\frac{p_3}{p_1}}+{\frac{p_3}{p_2}}}}\Bigl(\int\limits_{|t|<r}{\vert f(x-t) -
 f(x)\vert}^{p_1}dt\Bigr)^{\frac{p_3}{p_1}}\Bigl(\int\limits_{|t|<r}\vert
g(x-t)-g(x)\vert^{p_2}dt\Bigr)^{\frac{p_3}{p_2}} \nonumber 
\end{eqnarray}
\begin{eqnarray}
&&+{\frac {1}{r}}\Bigl(\int\limits_{|t|<r}{\vert
f(x)\vert}^{p_1}dt\Bigr)^{\frac{p_3}{p_1}}\Bigl(\int\limits_{|t|<r}\vert
g(x-t)-g(x)\vert^{p_2}dt\Bigr)^{\frac{p_3}{p_2}} \nonumber \\
& = & \frac {1}{r^{{\frac{p_3}{p_1}}}}\Bigl(\int\limits_{|t|<r}{\vert f(x-t) - f(x)\vert}^{p_1}dt\Bigr)^
{\frac{p_3}{p_1}}\frac {1}{r^{{\frac{p_3}{p_2}}}}\Bigl(\int\limits_{|t|<r}\vert
g(x-t)-g(x)\vert^{p_2}dt\Bigr)^{\frac{p_3}{p_2}} \nonumber \\
&&+{\frac {1}{r}}\Bigl[{\vert f(x)\vert}^{p_3}\Bigl(\int\limits_{|t|<r}1^{p_1}dt\Bigr)^{\frac{p_3}{p_1}}\Bigr]
\Bigl(\int\limits_{|t|<r}\vert g(x-t)-g(x)\vert^{p_2}dt\Bigr)^{\frac{p_3}{p_2}} \nonumber \\
&=&\frac {1}{r^{{\frac{p_3}{p_1}}}}\Bigl(\int\limits_{|t|<r}{\vert f(x-t) - f(x)\vert}^{p_1}dt\Bigr)^
{\frac{p_3}{p_1}}\frac {1}{r^{{\frac{p_3}{p_2}}}}\Bigl(\int\limits_{|t|<r}\vert
g(x-t)-g(x)\vert^{p_2}dt\Bigr)^{\frac{p_3}{p_2}} \nonumber \\
&&+2^{\frac {p_3}{p_1}}{\vert f(x)\vert}^{p_3}\frac {1}{r^{{\frac{p_3}{p_2}}}}
\Bigl(\int\limits_{|t|<r}\vert g(x-t)-g(x)\vert^{p_2}dt\Bigr)^{\frac{p_3}{p_2}} \nonumber
\end{eqnarray}
\begin{eqnarray}
&=&\Bigl[\Bigl(\frac {1}{r}\int\limits_{|t|<r}{\vert f(x-t) - f(x)\vert}^{p_1}dt\Bigr)^
{\frac{p_3}{p_1}}+2^{\frac{p_3}{p_1}}\vert
f(x)\vert^{p_3}\Bigr] \nonumber \\
&&\cdot \Bigl ({\frac {1}{r}}\int\limits_{|t|<r}\vert
g(x-t)-g(x)\vert^{p_2}dt\Bigr)^{\frac{p_3}{p_2}}\rightarrow 0
\;{\mbox  {as}  }\; r\rightarrow 0. \nonumber
\end{eqnarray}
Now we consider
\begin{eqnarray}
J&=&\frac {\vert g(x)\vert^{p_3}}{r}\int\limits_{|t|<r} \vert
f(x-t)-f(x)\vert^{p_3}dt. \nonumber
\end{eqnarray}
Since $p_3\leq p_1,$ it follows that $f\in L^{p_3}_{loc}$. By
Lemma 1, $J\rightarrow 0$ as $r\rightarrow 0$ and the Lemma is
proved because $K\rightarrow 0$ as $r\rightarrow 0$.
\begin{lemma} Let $f_i\in L^{p_i}_{loc}$, $i=1,...,n$, $\displaystyle\sum_{i=1}^n \frac {1}{p_i}=\frac
{1}{q_{n-1}}<1$, $F=\displaystyle\prod_{i=1}^n f_i$ and $x\in
\displaystyle\bigcap_{i=1}^n A_{f_i}^{p_i}$. Then $x\in
A^{q_{n-1}}_{F}$.
\end{lemma}
{\bf Proof:} Note that $F\in L^{q_{n-1}}_{loc}$. Put $\frac
{1}{p_1}+\frac {1}{p_2}=\frac {1}{q_1}$.  Lemma 3 implies that
$x\in A_{f_1f_2}^{q_1}$. Now we have $x\in A_{f_1f_2}^{q_1}$,
$x\in A_{f_i}^{p_i}$ for $i\geq 3$ and $\frac {1}{q_1}+\frac
{1}{p_3}+\cdots +\frac {1}{p_n}<1$. We put $\frac {1}{q_1}+\frac
{1}{p_3}=\frac {1}{q_2}$ and obtain that $x\in
A_{f_1f_2f_3}^{q_2}$. If we continue this procedure $n-1$ times at
the last iteration we have $\frac {1}{q_{n-2}}+\frac
{1}{p_n}=\frac {1}{q_{n-1}}$ and again, using Lemma 3 it follows
that $x\in A^{q_{n-1}}_{F}$.
\begin{lemma} Let $f_i\in L^{p_i}_{loc}$, $i=1,...,n$, $\displaystyle\sum_{i=1}^n \frac {1}{p_i}=1$,
and $F=\displaystyle\prod_{i=1}^n f_i$. Then
 $x\in \displaystyle\bigcap_{i=1}^n A_{f_i}^{p_i}$ implies that $x\in A_F$.
\end{lemma}
{\bf Proof:} The proof follows by Lemmas 4 and 2.

\section{An Approximation Lemma for BHT}

Our aim is to prove the following version of Lemma 1.2 Ch.VI in [8].
\begin{lemma} Let $f\in L^{p_1}(\mathbb{R})$, $g\in L^{p_2}(\mathbb{R})$, $p_1,p_2\in
[1,\infty),{\frac{1}{p_1}}+{\frac{1}{p_2}}=1$. Let $x\in
A^{p_1}_f\bigcap A^{p_2}_g$. Then $x\in A_{fg}$ and
$$\lim_{\varepsilon \rightarrow
0}\Bigl\{\int\limits_{-\infty}^{\infty}f(x-t)g(x+\alpha
t)\frac{t}{t^{2}+\varepsilon^{2}}dt-\int\limits_{0<\varepsilon\leq|t|}\frac{f(x-t)g(x+\alpha
t)}{t}dt\Bigr\}=0.$$
\end{lemma}

For the proof of Lemma 6 we need a version of Theorem 1.25, Ch. I
in [8]. This is Lemma 7:
\begin{lemma} Let $\varphi\in L^1(\mathbb{R})$, $\psi(x)=ess.\sup_{|t|\geq |x|}|\varphi(t)|, \;x\in
\mathbb{R}$. Assume that $\psi\in L^1(\mathbb{R})$. Let
$\varphi_\varepsilon=\frac{1}{\varepsilon}\varphi\bigr(\frac{\cdot}{\varepsilon}\bigl),\;\varepsilon>0.$
Assume $f\in L^{p_1}(\mathbb{R})$, $g\in L^{p_2}(\mathbb{R})$,
$p_1,p_2\in [1,\infty),{\frac{1}{p_1}}+{\frac{1}{p_2}}=1$, $x\in
A^{p_1}_f\bigcap A^{p_2}_g$. Then $x\in A_{fg}$ and
$$\lim_{\varepsilon \rightarrow
0}\int\limits_{\mathbb{R}}f(x-t)g(x+\alpha
t)\varphi_\varepsilon(t)dt=f(x)g(x)\int\limits_{\mathbb{R}}\varphi(t)dt.$$
\end{lemma}
{\bf Proof of Lemma 7:} First we prove the following assertion:

Let $\delta >0$. Then there exists $\eta>0$ such that
\begin{equation}\label{c}
\frac{1}{r}\int\limits_{|t|<r}|f(x-t)g(x+\alpha
t)-f(x)g(x)|dt<\delta\;{\mbox {if}}\;r\leq\eta.
\end{equation}

We have
$$
\frac{1}{r}\int\limits_{|t|<r}|f(x-t)g(x+\alpha t)-f(x)g(x+\alpha
t)+f(x)g(x+\alpha t)-f(x)g(x)|dt
$$
$$
\leq\frac{1}{r}\int\limits_{|t|<r}|f(x-t)-f(x)|\cdot|g(x+\alpha
t)|dt+\frac{1}{r}\int\limits_{|t|<r}|f(x)|\cdot|g(x+\alpha
t)-g(x)|dt=I_1+I_2.
$$
We estimate $I_1$ as follows:
$$
I_1\leq
\frac{1}{r}\Bigl(\int\limits_{|t|<r}|f(x-t)-f(x)|^{p_1}dt\Bigr)^{\frac{1}{p_1}}\Bigl(\int\limits_{|t|<r}|g(x+\alpha
t)-g(x)+g(x)|^{p_2}dt\Bigr)^\frac{1}{p_2}
$$
$$
\leq\Bigl(\frac{1}{r}\int\limits_{|t|<r}|f(x-t)-f(x)|^{p_1}dt\Bigr)^{\frac{1}{p_1}}\Bigl(\frac{1}{r}\int\limits_{|t|<r}|g(x+\alpha
t)-g(x)|^{p_2}dt+\frac{1}{r}\int\limits_{|t|<r}|g(x)|^{p_2}dt\Bigr)^\frac{1}{p_2}
.$$ Substituting $\alpha t=u$ we obtain
$$
I_1\leq\Bigl(\frac{1}{r}\int\limits_{|t|<r}|f(x-t)-f(x)|^{p_1}dt\Bigr)^{\frac{1}{p_1}}
\Bigl(\frac{\alpha}{\alpha
r}\int\limits_{|t|<\alpha r}|g(x-t)-g(x)|^{p_2}dt+C\Bigr)^\frac{1}{p_2}.
$$
Assumptions on $f$ and $g$ imply that the right hand side tends to zero as
$r\rightarrow 0.$ Similarly for $I_2$ we have the same inequality because
$$
I_2\leq
\frac{1}{r}\Bigl(\int\limits_{|t|<r}|f(x)|^{p_1}dt\Bigr)^{\frac{1}{p_1}}\Bigl(\int\limits_{|t|<r}|g(x+\alpha
t)-g(x)|^{p_2}dt\Bigr)^\frac{1}{p_2}
$$
$$
\leq C\Bigl(\frac{\alpha}{\alpha r}\int\limits_{|t|<\alpha
r}|g(x-t)-g(x)|^{p_2}dt\Bigr)^\frac{1}{p_2}.
$$
The last expression tends to zero as $r \rightarrow 0.$ Thus, $(\ref{c})$ is proved.

Let $a=\int\limits_{\mathbb{R}}
\varphi_\varepsilon(t)dt=\int\limits_{\mathbb{R}} \varphi(t)dt$.
We have
$$
\Bigl|\int\limits_{\mathbb{R}}f(x-t)g(x+\alpha
t)\varphi_{\varepsilon}(t)dt-af(x)g(x)\Bigr|\leq
$$
$$
\leq\Bigl|\int\limits_{|t|<\eta}\bigl[f(x-t)g(x+\alpha
t)-f(x)g(x)\bigr]\varphi_\varepsilon(t)dt\Bigr|
$$
$$
+\Bigl|\int\limits_{|t|\geq\eta}\bigl[f(x-t)g(x+\alpha
t)-f(x)g(x)\bigr]\varphi_\varepsilon(t)dt\Bigr|=I_1+I_2.
$$

 As in [8], p. 14, we will use the properties of the
function $\psi$. Function $\psi$ is radial ($\psi(x_1)= \psi(x_2)$
if $|x_1|=|x_2|$). Put $\psi_0(r)=\psi(x)$, $|x|=r$ then $\psi_0$
is a decreasing function of $|r|$.  We will use the fact that
\begin{equation}
\label{dva}
\lim_{t\rightarrow0}t\psi_0(t)\rightarrow 0 \mbox{ as } t\rightarrow 0
\mbox{ or } t\rightarrow\infty.
\end{equation}
Put
$$
d_+(t)=|f(x-t)g(x+\alpha t)- f(x)g(x)|,\quad t> 0,
$$
$$
d_-(t)=|f(x-t)g(x+\alpha t)- f(x)g(x)|,\quad t< 0.
$$

We have
$$d_{+}\leq d_{1,+}+d_{2,+},\;  d_{-}\leq d_{1,-}+d_{1,-},
$$
where
$$d_{1,+}=|f(x-t)g(x-t)- f(x)g(x)|,\; d_{2,+}=|f(x-t)g(x+\alpha t)- f(x-t)g(x-t)|,\; t>0,
$$
and $d_{1,-}$ and $ d_{2,-}$ are defined in the same way, for $t<0.$

Put
$$
D_{i,+}(t)=\int\limits_{0}^{t}d_{i,+}(s)ds,\quad t> 0, \quad{\mbox {and}}
\quad D_{i,-}(t)=\int\limits_{0}^{t}d_{i,-}(s)ds,\quad t< 0,\; i=1,2.
$$
We will estimate $D_{1,+}(t)$ and in the same way one can estimate $D_{1,-}(t)$.
For the estimates of $D_{2,+}(t)$ and $D_{2,-}(t)$ one has to make
the decomposition
$$d_{2,+}(t)\leq  |f(x-t)||g(x+\alpha t)-g(x)|,\;t>0,
$$
 similarly for
$d_{2,-}(t),\;t<0,
$
and to repeat the proof as for  $D_{1,+}(t)$, with a simple change of variable $\alpha t=-s$ in the integrals  related to $d_{2,+}$ and $d_{2,+}$.

We have
$$
I_1\leq\Bigl|\int\limits_{0}^{\eta}d_{1,+}(t)\varphi_\varepsilon(t)dt\Bigr|
+\int\limits_{0}^{\eta}d_{2,+}(t)\varphi_\varepsilon(t)dt\Bigr|
$$
$$
+\Bigl|\int\limits_{-\eta}^{0}d_{1,-}(t)\varphi_\varepsilon(t)dt\Bigr|
+\Bigl|\int\limits_{-\eta}^{0}d_{2,-}(t)\varphi_\varepsilon(t)dt\Bigr|
=I_{11}+I_{12}+I_{13}+I_{14}.
$$
We will estimate $I_{11}$;  the similar estimate holds  for $I_{12}, I_{13}$ and $I_{14}$.

For every $\varepsilon<1$,
$$
I_{11}=
\big|D_{1,+}(t)\frac{1}{\varepsilon}\psi_0(\frac{t}{\varepsilon})\Bigl|_0^\eta-
\int\limits_{0}^{\eta}D_{1,+}(t)d(\frac{1}{\varepsilon}\psi_0\bigl(\frac{t}{\varepsilon}\bigr))\big|
$$
$$\leq\big|\frac{D_{1,+}(\eta)}{\eta}\sup_{\varepsilon<1}\frac{\eta}{\varepsilon}\psi_0(\frac{\eta}{\varepsilon})-
\int\limits_{0}^{\eta/\varepsilon}D_{1,+}(\varepsilon t)\frac{1}{\varepsilon}
d(\psi_0(t))\big|.
$$
Note (cf. (\ref{dva})),
$$\frac{D_{1,+}(\eta)}{\eta}\rightarrow 0 \mbox{ as } \eta\rightarrow
0 \mbox{ and }
\sup_{\varepsilon<1}\frac{\eta}{\varepsilon}\psi_0(\frac{\eta}{\varepsilon})<\infty.$$
Now, we estimate the second integral. By the assumption,
$$D_{1,+}(\varepsilon t) \leq\varepsilon \eta, t<\eta.$$
With this we have:
$$|\int\limits_{0}^{\eta}D_{1,+}(t)d\Bigl(\frac{1}{\varepsilon}\psi_0\bigl(\frac{t}{\varepsilon}\bigr)\Bigr)|\leq C \eta
\int\limits_{0}^{\infty}|d\psi_0(s)|,
$$
Thus, for every $\delta>0$, there exists $\eta_0$ such that
$$I_{11}<\delta, \; \eta<\eta_0.
$$
The same holds for $I_{12}, I_{13}, I_{14}.$
Let us estimate $I_{2}=I_{21}+I_{22}+I_{23}+I_{24}$. Actually, we will estimate only $I_{21}$ since the
other parts can be estimated in the same way.
$$
I_{21}\leq
\sup_{\varepsilon<1,t>\eta}\frac{1}{t}\frac{t}{\varepsilon}|\varphi(\frac{t}{\varepsilon})|
\int\limits_{|t|>\eta}|f(x-t)g(x- t)-f(x)g(x)|dt $$
$$\leq
\frac{1}{\eta}
\sup_{\varepsilon<1,t>\eta}|\frac{t}{\varepsilon}\varphi(\frac{t}{\varepsilon})|
\int\limits_{|t|>\eta}|f(x-t)g(x-t)-f(x)g(x)|dt.
$$
Now the assertion follows by (\ref{c}).

 Thus we have that for every $\delta>0$ there exists $\varepsilon_0$
 such that
$$\Bigl|\int\limits_{\mathbb{R}}f(x-t)g(x+\alpha
t)\varphi_\varepsilon(t)dt-af(x)g(x)\Bigr|<C\delta.
$$
This completes the proof of Lemma 7.

{\bf Proof of Lemma 6:} The
proof follows by the use of the same idea as in Lemma 1.2 Ch.VI in
[8]. Put
$$\varphi(t)=\left\{\begin{array}{rl}
\frac{t}{t^2+1}-\frac{1}{t}, & |t|\geq 1 \\
\frac{t}{t^2+1}, & |t|< 1,
\end{array}
\right.$$
$\varphi_\varepsilon(t)=\varepsilon^{-1}\varphi\bigr(\frac{t}{\varepsilon}\bigl),\;\varepsilon>0,\;t\in
\mathbb{R}$ and $\psi(x)=\sup_{|t|\geq |x|}|\varphi(t)|, \;x\in \mathbb{R}$. Note
$\psi\in L^1(\mathbb{R})$ and
$$\int\limits_{\mathbb{R}}f(x-t)g(x+\alpha t)
\frac{t}{t^{2}+\varepsilon^{2}}dt-
\int\limits_{0<\varepsilon\leq|t|}\frac{f(x-t)g(x+\alpha t)}{t}dt$$
$$=\int\limits_{\mathbb{R}}f(x-t)g(x+\alpha t)\varphi_\varepsilon(t)dt.$$
Now the proof follows by Lemma 7 and the observation (as in [8],
p.218) that
$$\lim_{\varepsilon \rightarrow
0}\int\limits_{\mathbb{R}}f(x-t)g(x+\alpha
t)\varphi_\varepsilon(t)dt=f(x)g(x)\int\limits_{\mathbb{R}}\varphi(t)dt=0,$$ since
$\int\limits_{\mathbb{R}}\varphi(t)dt=0$.

{\bf Remark.}
With the same arguments as for Lemma 6 we have

{\it
Let $f\in L^{p_1}$, $g\in L^{p_2}$, where
${\frac{1}{p_1}}+{\frac{1}{p_2}}={\frac{1}{p_3}}<1$, $x\in
A^{p_1}_f\bigcap A^{p_2}_g$. Then it follows that $x\in
A^{p_3}_{fg}$ and
$$
\lim_{\varepsilon \rightarrow 0}\Bigl\{\int\limits_{\mathbb{R}}f(x-t)g(x+\alpha
t)\frac{t}{t^{2}+\varepsilon^{2}}dt-\int\limits_{0<\varepsilon\leq|t|}\frac{f(x-t)g(x+\alpha
t)}{t}dt\Bigr\}=0.
$$
}

{\bf Proof of  Theorem 1.1}. Let us return to the  proof of Theorem 1.1. in Section 1.

By
Lemma 6, it follows that $(\ref{a})$ tends to zero as $\varepsilon\rightarrow 0.$
Since
$$
\lim_{\varepsilon\rightarrow 0}\int\limits_{-\infty}^{\infty}[f(x-t)g(x+\alpha
t)-f(x)g(x)]\frac {\varepsilon dt}{t^2+\varepsilon^2}
$$
$$
=
\lim_{\varepsilon\rightarrow 0}\int\limits_{-\infty}^{\infty}[f(x-\varepsilon t)g(x+\alpha\varepsilon
t)-f(x)g(x)]\frac { dt}{t^2+1}
$$
by the Lebesque theorem,  $(\ref{b})$ tends to zero as $\varepsilon \rightarrow 0.$

This completes the proof of Theorem 1.1.

\section{Appendix}
As a consequence of Theorem 1.1, we recall and reprove parts (iii)
and (iv) of Theorem 5 in [1].

\begin{theorem}
Let $f\in {\cal D}'_{L^{2}}(\mathbb{R})$, $g\in {\cal D}_{L^{\infty}}(\mathbb{R})$, respectively $f\in
{\cal D}'_{L^{q}}(\mathbb{R})$, $g\in {\cal D}_{L^{p_2}}(\mathbb{R})$, $1<p_1,p_2,\;{p=\frac{p_1p_2}{p_1+p_2}},\;{q=\frac{p}{p-1}}$ and ${q_1=\frac{p_1}{p_1-1}}$. Then in the sense of convergence in ${\cal D}'_{L^2}$, respectively ${\cal D}'_{L^{q_1}}$ holds:
$$
f(x)g(x)=\frac{i}{\pi}(\lim_{\varepsilon\rightarrow
0}H^*_{\alpha,\varepsilon}(f,g)(x)-H^*_{\alpha}(f,g)(x)),
$$
where
$$
<H^*_{\alpha,\varepsilon}(f,g),\psi>=<f,H_{\alpha,\varepsilon}(\psi,g)>,\quad \psi \in {\cal D}_{L^2}, {\mbox {respectively}}
\quad \psi \in {\cal D}_{L^{p_1}}
$$
and
$$
<H^*_{\alpha}(f,g),\psi>=<f,H_{\alpha,g}\psi>,\quad \psi \in {\cal D}_{L^2}, {\mbox {respectively}}
\quad \psi \in {\cal D}_{L^{p_1}}
$$
are defined in [1], Definition 6 and 1.
\end{theorem}

{\bf Proof:} We use the simple property
$$
\Bigl( \frac{d}{dx} \Bigr)^m H_{\alpha,\varepsilon}(f,g)(x)=\sum_{k=0}^m{m\choose k}(H_{\alpha,\varepsilon}(f^{(k)},g^{(m-k)})(x))
$$
and
\begin{eqnarray}
\lim_{\varepsilon \rightarrow 0^+}<H^*_{\alpha,\varepsilon}(f,g),\psi>&=&\lim_{\varepsilon \rightarrow 0^+}<f,H_{\alpha,\varepsilon}(g,\psi)>\nonumber \\
&=&<f,H_{\alpha}(g,\psi)-i\pi g\psi> \nonumber \\
&=&<f,H_{\alpha}(g,\psi)>-<f,i\pi g\psi> \nonumber \\
&=&<H^*_{\alpha}(f,g),\psi>-<i\pi fg,\psi>\nonumber \\
&=&<H^*_{\alpha}(f,g)-i\pi fg,\psi>. \nonumber
\end{eqnarray}
which is the consequence of Lemma 6 for $g$ and $\psi$. This implies the above assertions.

\end{document}